\documentclass[11pt,a4paper]{article}
\usepackage[utf8]{inputenc}
\usepackage[spanish,english]{babel}
\usepackage{amsmath,amssymb,amsthm}
\usepackage{accents}

\usepackage{graphicx}
\usepackage{hyperref}
\usepackage{geometry}
\title{\textbf{\Huge{Khajuraho's magic square \\
is an hypercube}\footnote{This is a translation to English of the article 
{\em Le carré magique de Khajuraho est un hypercube}, published in 2020 in the CNRS disemination journal 
{\em Images des Math\'ematiques} \cite{Na20}, and translated to Spanish in the mirror site {\em Paisajes Matemáticos} 
\cite{Na'20}. I strongly thank the participants of the special session on history of mathematics of the ICM 2026, particularly 
Clemency Montelle, for pushing me to make this translation available.}}}
\author{Andrés Navas\footnote{Except for the first two pictures (widely available on the internet), all the figures of this 
article were produced by María José Moreno, to whom I would like to extend my warm gratitude.}}
\date{}

\begin{document}

\maketitle

\begin{figure}[htbp]
    \centering
    \includegraphics[width=0.6\textwidth]{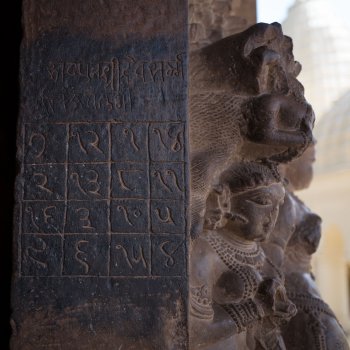}
\end{figure}

\begin{center}
The panmagic square engraved in a temple in Khajuraho, India, and studied by $\mathrm{N\overline{a}r\overline{a}ya\underaccent{.}na \,\, 
Pa\underaccent{.}n\underaccent{.}dita}$ in the 14th century, closely guards a secret: its group of symmetries is isomorphic to that of the 
hypercube, the four-dimensional analogue of the standard cube.
\end{center}

\vspace{0.4cm}

\noindent{\bf \Large A Formidable Ancestral Object}

\vspace{0.3cm}

The town of Khajuraho in central India, declared a UNESCO World Heritage Site in 1986, is famous for its temples decorated with motifs blending divine 
and erotic (or, rather, tantric) elements. In the porch of the Parshwanath temple stands a striking engraving dating back to at least the 12th century: 
a number grid including the entries  $1, 2, \dots, 16$ (written in Sanskrit\footnote{After the publication of the article in 2019, I learned from 
Gautami Bhowmik that these are not Sanskrit characters but coming from a close language.}).

\newpage

\begin{figure}[htbp]
    \centering
    \includegraphics[width=0.4\textwidth]{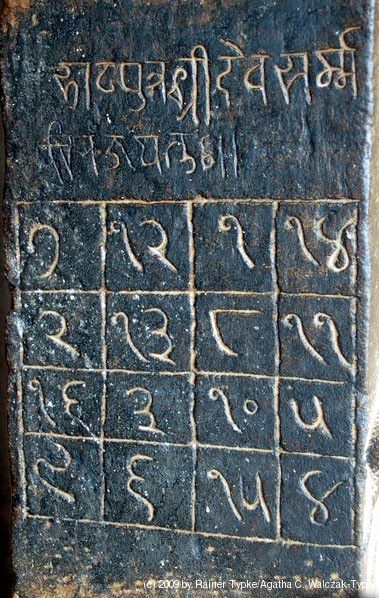}
\end{figure}

These appear arranged in such a way that certain ``(pan)magic'' properties are satisfied: in each illustration, the numbers in the cells of the same color sum up to exactly $34$.

\begin{figure}[htbp]
    \centering
    \includegraphics[width=0.9\textwidth]{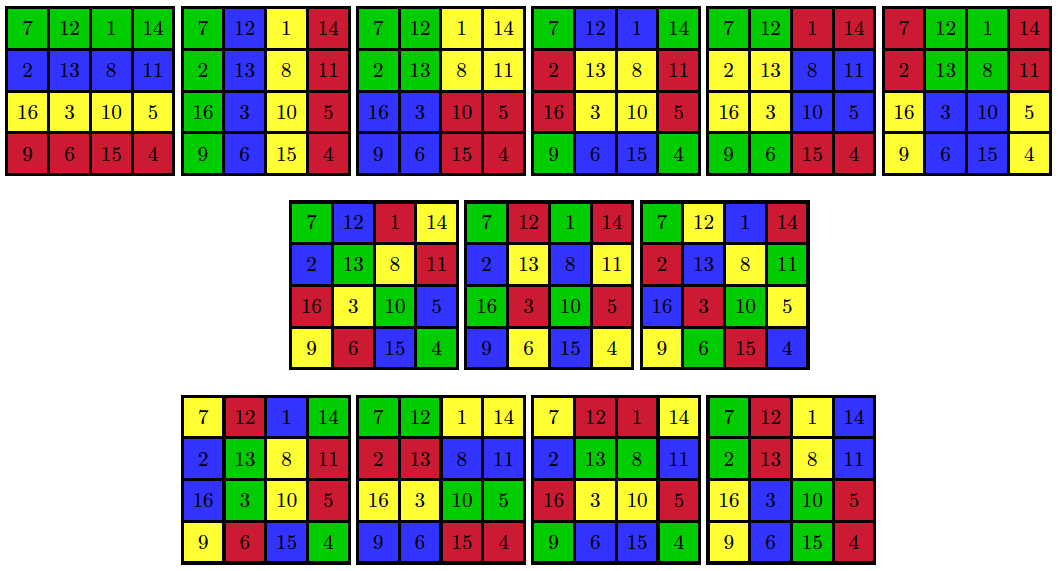}
\end{figure}

This grid, commonly called the \emph{Chautisa Yantra}, has already been the subject of two articles in \emph{Images des Math\'ematiques}. On the one hand, Gautami Bhowmik 
addresses the 14th-century study by the celebrated mathematician $\mathrm{N\overline{a}r\overline{a}ya\underaccent{.}na \,\, Pa\underaccent{.}n\underaccent{.}dita}$  
concerning its combinatorial properties 
and its relationship with the movement of a chess knight \cite{Bo18}. On the other hand, I analyzed the most general possible configurations of this type (with non-necessarily 
consecutive numbers) via elementary linear algebra \cite{Na19}. My objective now is to introduce another element for the discussion and understanding of this marvelous 
object: its group of symmetries. This is inspired by the following universal principle:

\begin{figure}[htbp]
    \centering
    \includegraphics[width=0.55\textwidth]{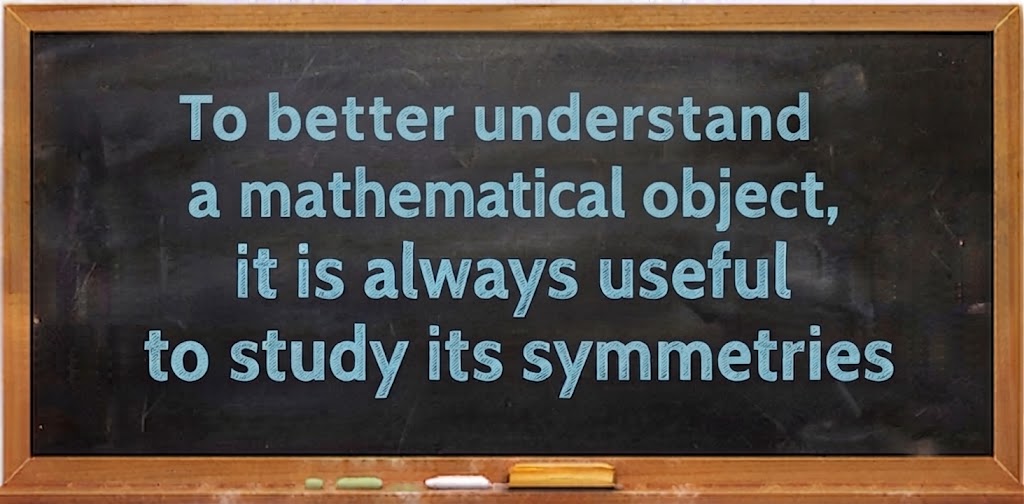}
\end{figure}

\newpage

\begin{quote}
\textbf{A Brief Reminder on Groups}

Recall that a \emph{bijective transformation} (or simply a \emph{bijection}) of a space is a transformation that sends each ``point'' 
of the space to another point of the same space in such a way that:

\vspace{0.1cm}   

--  different points are sent to different points (\emph{injectivity});

\vspace{0.1cm}   

-- every point in the space is the result of the transformation of another point (\emph{surjectivity}).

\vspace{0.1cm}

When the space is finite, the first condition implies the second, but this ceases to hold for infinite spaces.

Bijections can be ``composed'' (or ``multiplied''): if $S$ and $T$ are two bijections, then their ``product'' $S \circ T$ transforms 
a point $x$ of the space into $S(T(x))$, which corresponds to the result of the transformation $S$ applied to the point resulting 
from applying the transformation $T$ to the point $x$.

A \textbf{\em transformation group} is a set of bijections such that:

\vspace{0.1cm}   

-- the composition of two elements of the group still belongs to the group;
 
 \vspace{0.1cm}    
 
-- the \emph{identity} (that is, the transformation leaving all points intact) belongs to the group;
  
 \vspace{0.1cm}    
 
-- for each transformation $T$ in the group, its \emph{inverse} $T^{-1}$ also belongs to it (recall that the 
 inverse of $T$ is the transformation that returns each point to its original position before $T$ acts).

\vspace{0.1cm}

Note that composition is necessarily \emph{associative}, that is,
$$R \circ (S \circ T) = (R \circ S) \circ T.$$

The group of all bijective transformations is called the \emph{permutation group} of the set. In many cases, the transformations worked with have a geometric origin. This is why one colloquially speaks of the ``symmetry group'' of an object (space).

The notion of a transformation group leads to that of an \textbf{\em abstract group}. An abstract group is a set of elements for which there is an ``internal product rule'' $ab$ satisfying the same formal properties above:
\begin{itemize}
    \item Associativity: $a(bc) = (ab)c$;
    \item Neutral element: there exists an element (denoted $e$) such that $ae = ea = a$ for every element $a$ of the group;
    \item Inverse element: for every element $a$ of the group there exists an element (denoted $a^{-1}$) such that $aa^{-1} = a^{-1}a = e$.
\end{itemize}

A fundamental theorem (both profound and elementary) due to Arthur Cayley states that every abstract group can be realized 
as a transformation group; in fact, it is a transformation group of itself.

The \emph{structure of a group} refers to the group as an abstract object. In this sense, two groups are equivalent (or \emph{isomorphic}) if there exists a correspondence between their elements that preserves the properties related to multiplication in each of them (for example, the correspondence must take the product of two elements to the product of their corresponding elements...).
\end{quote}

\section*{An Associated Group?}

Recall from the previous article \cite{Na19}
that a $4 \times 4$ panmagic square is a square filled with $16$ numbers such that the sums along the rows, 
columns, and diagonals (including broken diagonals) are all equal. The ``general formula'' for a square satisfying these properties is the following:

\begin{figure}[htbp]
    \centering
    \includegraphics[width=0.55\textwidth]{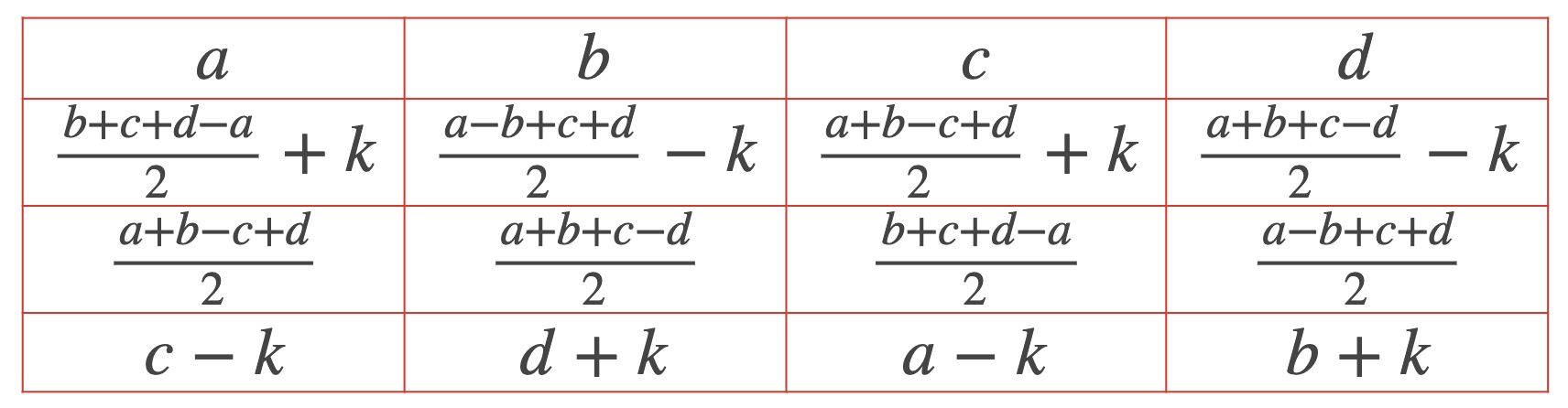}
\end{figure}

It is quickly verified that in a square of this type, the $52$ colored combinations mentioned above all yield the same sum (whose value is $a+b+c+d$). In other words, the $16$ initial equalities ($4$ rows, $4$ columns, $8$ diagonals) imply another $36$.

Consider the permutations of the cells that transform any panmagic square into another panmagic square. For example, this is the case for a $90^\circ$ rotation, since rows transform into columns and vice versa, while descending diagonals swap with ascending ones. Since all associated sums were originally equal, they remain so after the movement. Thus, a panmagic square subjected to this movement gives rise to a new panmagic square.

Panmagic transformations obviously form a group: the \textbf{\em panmagic group}. How many of these are there? What structure does this group have? We will answer these questions after a detailed analysis.

In addition to rotations, it is easy to imagine other ``panmagic'' transformations, such as reflections with respect to the vertical and horizontal axes, or with respect to the diagonals. But these are not the only ones. For example, one can observe from the formulas above that the permutation illustrated below is also panmagic:

\begin{figure}[htbp]
    \centering
    \includegraphics[width=0.45\textwidth]{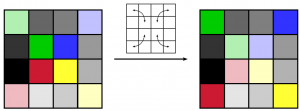}
\end{figure}

This is an \emph{involution}, that is, a permutation of order $2$ (meaning that applying it twice returns everything to its original position). To be even more surprised, 
here are illustrations of four other rather less evident panmagic transformations. For them, the orders are respectively $3$, $4$, $6$, and $8$.\footnote{The algebraic 
description of the panmagic group further described implies that the orders of its non-trivial elements can be equal to  $2, 3, 4, 6,$ and $8$.}  (Note that the element of 
order $4$ depicted below is a cyclic permutation of rows.) A nice exercise is to apply them to the Chautisa Yantra to obtain other squares with the same panmagic properties.

\newpage 

\begin{figure}[htbp]
    \centering
    \includegraphics[width=0.485\textwidth]{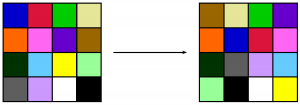}\\
    \vspace{0.1cm}
    \includegraphics[width=0.485\textwidth]{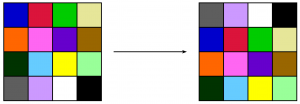}\\
    \vspace{0.1cm}
    \includegraphics[width=0.485\textwidth]{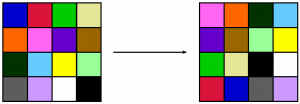}\\
    \vspace{0.1cm}
    \hspace{0.02cm} \includegraphics[width=0.485\textwidth]{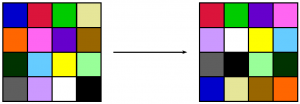}
\end{figure}

\section*{The Symmetries of the Lo Shu Square}

Determining the structure of the panmagic group is not easy. To develop an attack strategy, let us consider a simpler case: 
that of $3 \times 3$ magic squares. These are grids filled with numbers such that the sums along the rows, columns, and 
two diagonals are equal (this time, we do not consider broken diagonals). Again, we copy the general formula:

\begin{figure}[htbp]
    \centering
    \includegraphics[width=0.5\textwidth]{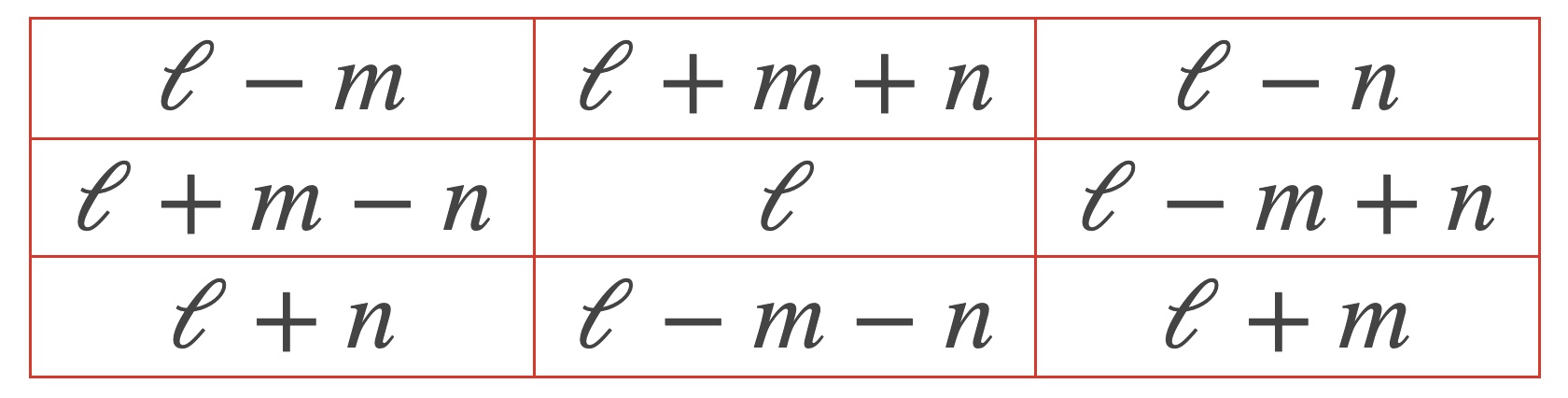}
\end{figure}

Clearly, the eight symmetries alluded to earlier (four rotations —one of them trivial— and four reflections) transform one magic configuration into another. 
Are there any other permutations with this property? The negative answer holds. Indeed, any transformation of this type must take the square below 
(called the \emph{Lo Shu} square) into another magic square with entries $1, 2, \dots, 9$.

\begin{figure}[htbp]
    \centering
    \includegraphics[width=0.12\textwidth]{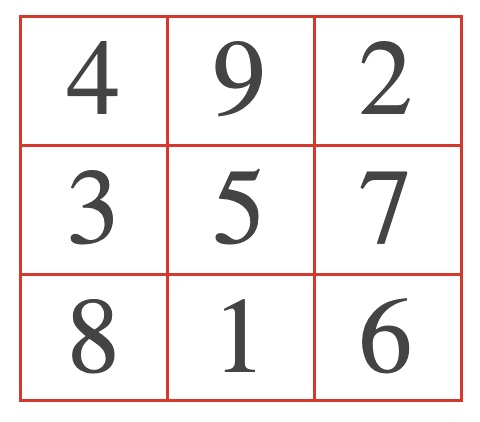}
\end{figure}

However, it is quickly verified with the formulas above that there are only $8$ squares with these characteristics. They are the ones obtained 
with the corresponding assignments of $m$ and $n$ (note that, necessarily, $\ell = 5$).

\newpage 

\begin{figure}[htbp]
    \centering
    \includegraphics[width=1\textwidth]{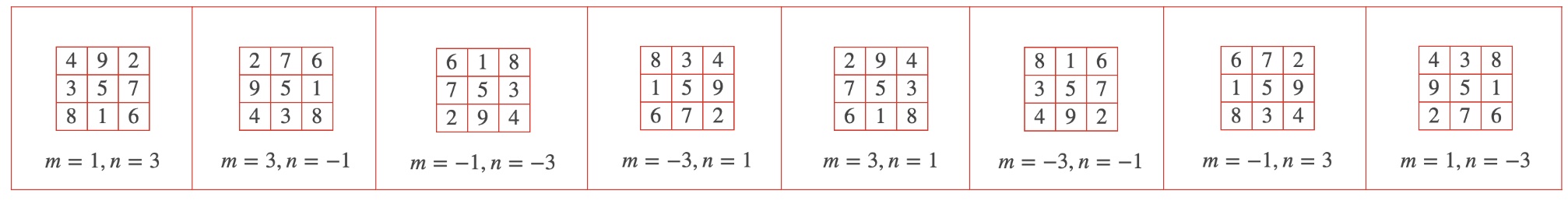}
\end{figure}

Therefore, the ``magic group'' (for $3 \times 3$ grids) consists of $8$ elements. Surely you have already recognized it: this is the famous dihedral group $D_4$, 
which coincides with the symmetry group of a square.

\section*{The Symmetries of the Chautisa Yantra}

To calculate the panmagic group, we will proceed analogously. First, we will seek to enumerate all panmagic squares that use the entries $1, 2, \dots, 16$, and then construct 
enough panmagic transformations to obtain all of them starting from the Chautisa Yantra. As we will see, rethinking the geometry of the ``square'' will prove 
fundamental. The final conclusion will be the following: just as the magic group is isomorphic to the symmetries of the square, \textbf{the panmagic group is isomorphic to the 
symmetry group of the hypercube.}

Let's continue~!

\subsection*{A New Geometry}

To begin, let us observe the combinations involved in the sums that take the same value in any panmagic square:

\begin{figure}[htbp]
    \centering
    \includegraphics[width=0.59\textwidth]{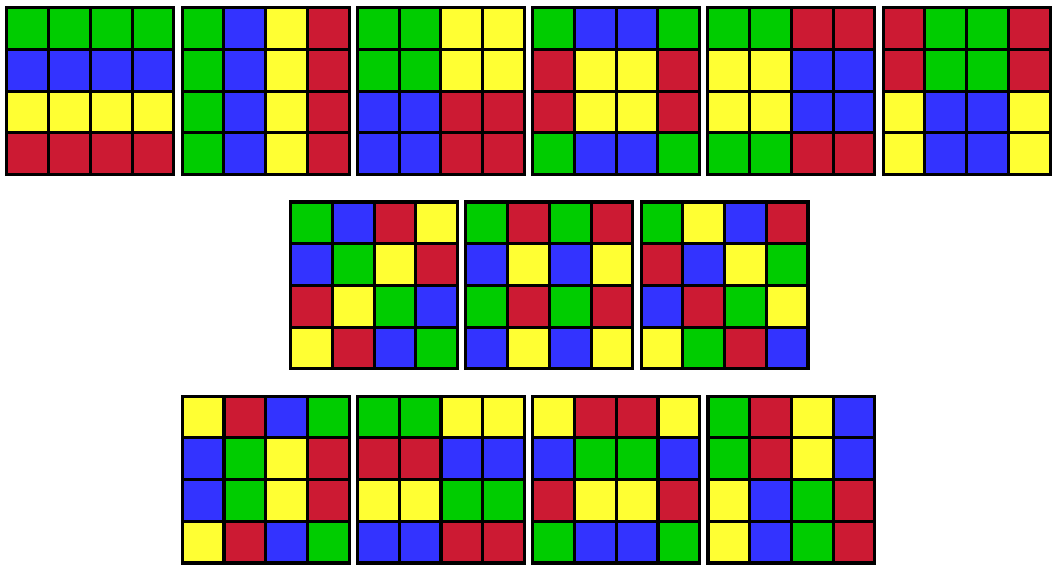}
\end{figure}

Is there any common pattern in them? The answer is a resounding YES, but to detect it, it is necessary to reinterpret things geometrically. The key is to think of a natural 
geometry for the $4 \times 4$ grid: the one in which the distance between two positions is the minimum number of cell boundaries crossed when traversing a path from 
one to another using only horizontal or vertical segments. Warning: when a horizontal path exits to the right or left, we think of it as entering on the other side at the 
height of the corresponding cell; likewise, if it exits downwards or upwards, it reappears on the opposite side.

The relevance of this new geometry associated with ``path distance'' (referred to simply as ``Distance'' in what follows) is made clear in the following statement.

\vspace{0.25cm}

\noindent{\bf Claim.}
\emph{Every cell permutation that preserves Distance is panmagic.}

\vspace{0.25cm}
The proof is given below. To get a feel for this, it is a good exercise to verify that the five permutations illustrated above preserve the Distance —that is, they 
map two squares at a given Distance to two squares located at the same Distance. In contrast, the permutation illustrated below does not preserve the Distance.

\begin{figure}[htbp]
    \centering
    \includegraphics[width=0.5\textwidth]{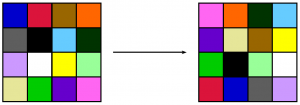}
\end{figure}

\noindent{\bf Proof of the Claim.} 
It suffices to analyze the $52$ sums appearing in a panmagic square in terms of Distance. Indeed, the corresponding quadruples can be classified into three types:
\begin{itemize}
    \item Type $(1,1,2)$: starting at any position, the other involved cells are at Distance $1$, $1$, and $2$ from it. There are $24$ configurations of this type.
    \item Type $(2,2,4)$: starting at any position, the other cells are at Distance $2$, $2$, and $4$ from it. There are $12$ configurations of this type.
    \item Type $(1,3,4)$: starting at any position, the other cells are at Distance $1$, $3$, and $4$ from it. There are $16$ configurations of this type.
\end{itemize}
This concludes the enumeration of the $52$ sums. If a permutation preserves Distance, it will take a configuration of one type to another of the same type. 
Since the sums of the entries in the original configurations were the same, they will also be in the new configurations. Consequently, the square remains panmagic.

\subsection*{The 384 Different Versions of the Chautisa Yantra}

$\mathrm{N\overline{a}r\overline{a}ya\underaccent{.}na \,\, Pa\underaccent{.}n\underaccent{.}dita}$ proved in his book $\mathrm{Ga\underaccent{.}nitakaumud\overline{\i}}$ 
that there are exactly $384$ panmagic squares involving each of the numbers $1, 2, \dots, 16$ (see \cite{Bo18} on this). 
We will see below that the Distance introduced above allows us to give another proof of this theorem, much simpler, 
based on the general formula obtained above.

\vspace{0.35cm}

\noindent{\bf Theorem
[$\mathrm{N\overline{a}r\overline{a}ya\underaccent{.}na \,\, Pa\underaccent{.}n\underaccent{.}dita}$].}
\emph{There exist $384$ panmagic squares using each of the entries $1, 2, \dots, 16$.}

\vspace{0.35cm}

Before giving the proof, we highlight two fundamental aspects arising from it.

\vspace{0.2cm}

\noindent 
\textbf{An algorithm to construct all panmagic squares:} Place $1$ in any position. Then, choose three of the four positions 
at Distance $1$ from it, and place the numbers $15, 14$, and $12$ there. Finally, complete the square so that it becomes 
panmagic using each number $1, 2, \dots, 16$ exactly once (there is a unique way of doing this).

\vspace{0.25cm}

\noindent 
\textbf{Counting from the algorithm:} Assuming the validity of the previous algorithm, to place $1$ we have $16$ different cells. 
For $15, 14, 12$, there are $4 \times 3 \times 2 = 24$ combinations. Thus, there are $16 \times 24 = 384$ panmagic squares.

\vspace{0.3cm}

\noindent{\bf Proof of $\mathrm{N\overline{a}r\overline{a}ya\underaccent{.}na \,\, Pa\underaccent{.}n\underaccent{.}dita}$'s theorem.} 
The maximum distance between two boxes is equal to $4$. In fact, given any square, there is a single square at distance exactly $4$ 
from it. We will say that two such boxes are {{opposite}}. In the picture below, opposite boxes have the same color. 

\begin{figure}[htbp]
    \centering
    \includegraphics[width=0.15\textwidth]{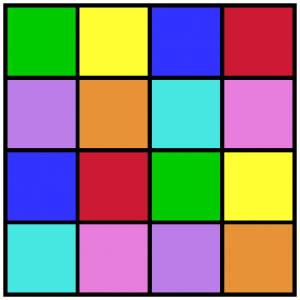}
\end{figure}

\newpage

Let us now notice that 
{{in every panmagic square, the sum of the entries in opposite squares is the same}} (equal to half of the ``magic number'' $a + b + c + d$). 
This can be verified by simple inspection from the general formula previously given. In particular, if we use (all) the numbers $1, 2, \ldots, 16$, 
then this sum is equal to $\frac{34}{2} = 17$.

\vspace{0.15cm}

Looking at the colored configurations above, we note that:

\vspace{0.15cm}

\noindent -- if two squares are at Distance $1$, then they appear linked in four of the $52$ sums; 

\vspace{0.15cm}

\noindent -- if two squares are at Distance $2$, then they are linked in two of these sums; 

\vspace{0.15cm}

\noindent -- if two squares are at Distance $3$, then they are linked in only a single sum;

\vspace{0.15cm}

\noindent -- if two squares are at Distance $4$, then they are linked in five sums.

\vspace{0.15cm}

With this information, we will prove that in a panmagic square that uses each number $1, 2, \ldots, 16$ once, the entries $15$, 
$14$, and $12$ must be placed in squares at Distance $3$ from the one where $16$ is located, and thus at Distance $1$ from the one 
where $1$ is located.

To verify this, note that there is only one way to complete the sum $16 + 15 = 31$ with two positive numbers to obtain $34$, 
namely $1 + 2$. But if $16$ and $15$ were at a Distance different from $3$, then they would have to be completed in at least 
two different ways... A similar argument applies for $14$, since $16 + 14 = 30$ can only be completed with $1 + 3$ (the complement 
$2 + 2$ is not valid because distinct additions are required). Unfortunately, for $12$, the argument is more elaborate (and a bit tedious), 
since $16 + 12 = 28$ can be completed in exactly two ways: $4 + 2$ and $1 + 5$.

If $16$ and $12$ were not at Distance $3$, then their 
Distance would have to be equal to $2$ (for a Distance of $1$ or $4$, they would have to be completed in more different ways, which is 
not the case). Consider the case below, all others being analogous (in fact, one can quickly reduce the general case to this one by using 
the action of permutations that preserve Distance, which we have already introduced).

\begin{figure}[htbp]
    \centering
    \includegraphics[width=0.175\textwidth]{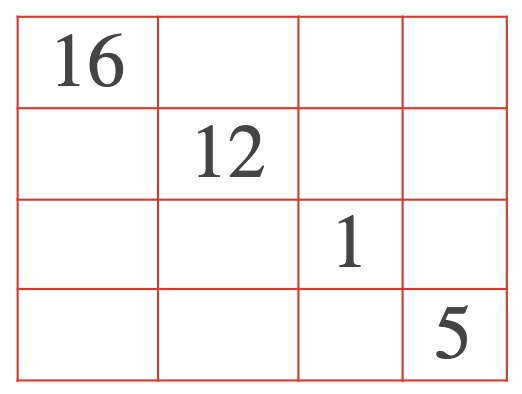}
\end{figure}

The entries $2$ and $4$ must appear below and 
to the right of $16$ (in order to obtain a sum of $34$ between the four top-right squares). Since both cases are analogous, we will assume 
that $4$ is below $16$, and $2$ is to the right, as illustrated below (again, it is possible to move from one case to the other by a reflection 
across the diagonal).

\begin{figure}[htbp]
    \centering
    \includegraphics[width=0.175\textwidth]{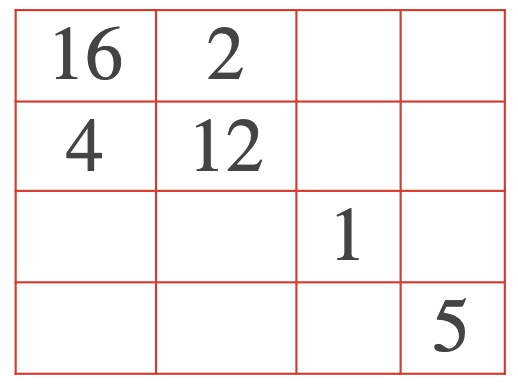}
\end{figure}

This forces $13$ and $15$ to appear in the corresponding opposite positions, thus ``semi-completing" the required 
sum of $17$:

\begin{figure}[htbp]
    \centering
    \includegraphics[width=0.175\textwidth]{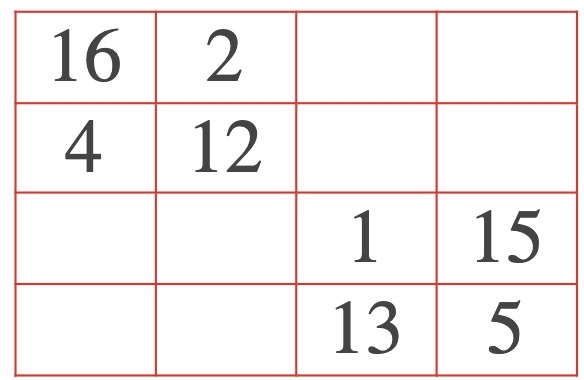}
\end{figure}

\newpage

Since we know that $14$ must appear in a position at Distance $3$ from $16$, the two cases below are the only possible 
ones. However, in both cases, the row that already has three complete entries must be filled with a $4$, and $4$ has already been used...

\begin{figure}[htbp]
    \centering
    \includegraphics[width=0.45\textwidth]{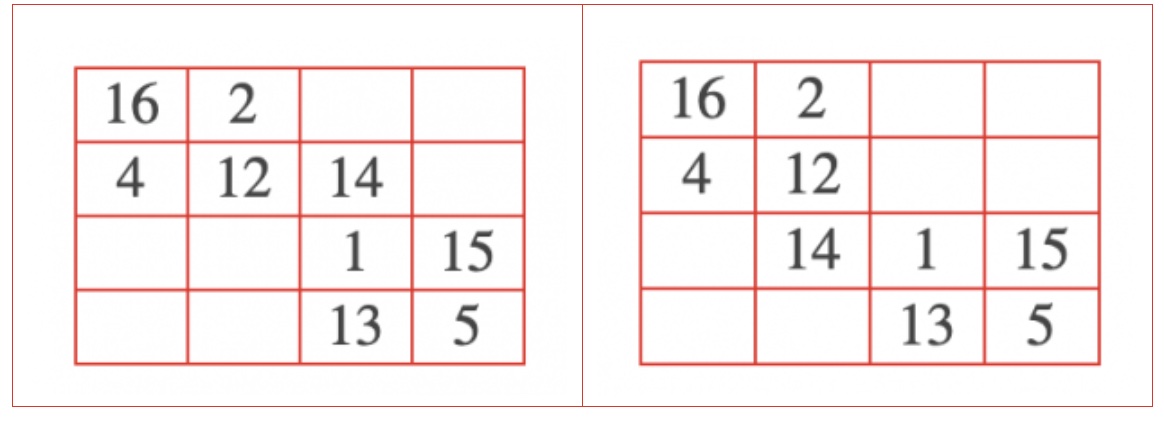}
\end{figure}

To conclude, one quickly verifies that for each configuration in which $15$, $14$, and $12$ appear around $1$ (and $16$ in a position 
opposite to it) there exists a unique way to complete the square with the entries $1, 2, \ldots, 16$ while respecting the panmagic 
properties (in fact, it is not necessary to check all possibilities, since appropriate Distance-preserving permutations can be used 
to simplify). An example is illustrated below.

\begin{figure}[htbp]
    \centering
    \includegraphics[width=0.65\textwidth]{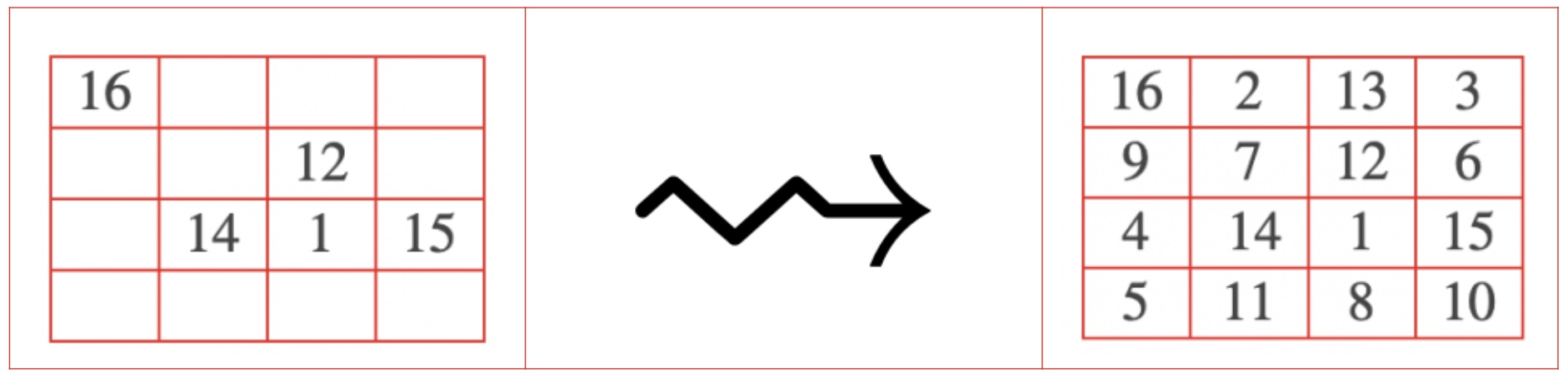}
\end{figure}

\section*{Finally, the Hypercube~!}

In our discussion, the geometry inherited from the Distance was essential. Is there a more concrete way to visualize the underlying mathematical object? 
Certainly yes: consider the illustration below, where each cell is indexed by four numbers equal to $0$ or $1$ each. Note that the indices from one row to 
another are the same except for one entry, which systematically changes. The same change occurs from one column to an adjacent one. In this way, the 
Distance between two cells is nothing more than the number of different digits between the indices of these two cells. For example, the indices of two 
opposite cells are completely different from one another: all four digits have changed.
\begin{figure}[htbp]
    \centering
    \includegraphics[width=0.32\textwidth]{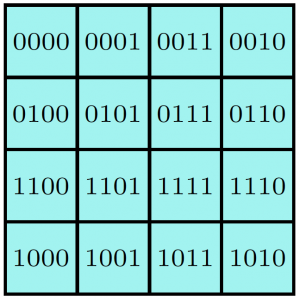}
\end{figure}

The object we are considering is then identified with the set of indices, which corresponds to the product $\{ 0,1 \}^4$. 
Even though it is still difficult to visualize this object, by replacing the exponent $4$ with a smaller number, everything becomes a little clearer:

\newpage

\begin{figure}[htbp]
    \centering
    \includegraphics[width=0.75\textwidth]{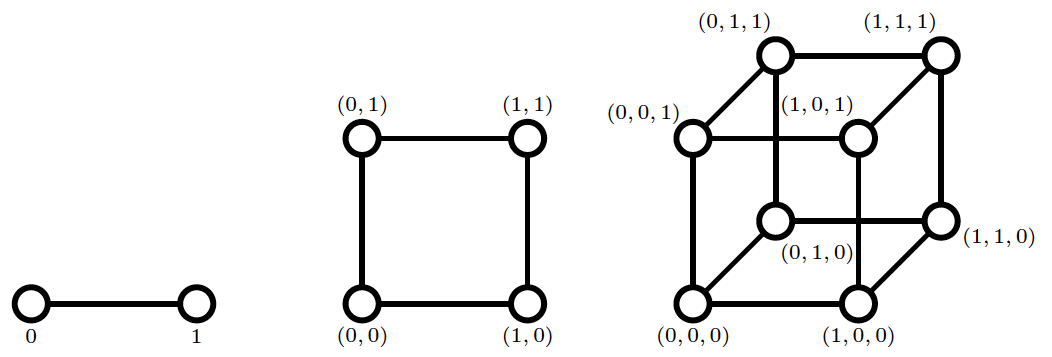}
\end{figure}

The set $\{ 0,1 \}^4$ is nothing more than the set of vertices of the object that must follow in the list 
beginning with the interval, the square, and the cube. This object is called the {{hypercube}} (or tesseract). 
Although it is an object of dimension $4$, we can ``project'' it into our $3$-dimensional space and then visualize it on a 
$2$-dimensional screen as follows:  

\begin{figure}[htbp]
    \centering
    \includegraphics[width=0.2\textwidth]{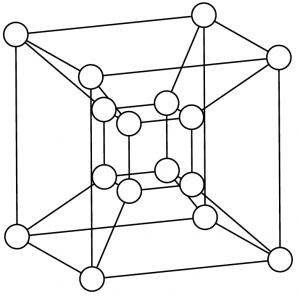}
\end{figure}

Observe that each vertex is connected by an edge to four other vertices 
(its ``neighbors''), in the same way that on the $4 \times 4$ square, each cell is at Distance $1$ from four other cells.  

Finally, with a little more familiarity with this structure, one will recognize that the list of $13$ hypercubes with colored vertices 
below corresponds to the list of $13$ colored $4 \times 4$ squares at the top (the configurations of type $(1,1,2)$ are colored 
in the same color in the first six, those of type $(2,2,4)$ in the middle three, and those of type $(1,3,4)$ in the last four).

\begin{figure}[htbp]
    \centering
    \includegraphics[width=0.99\textwidth]{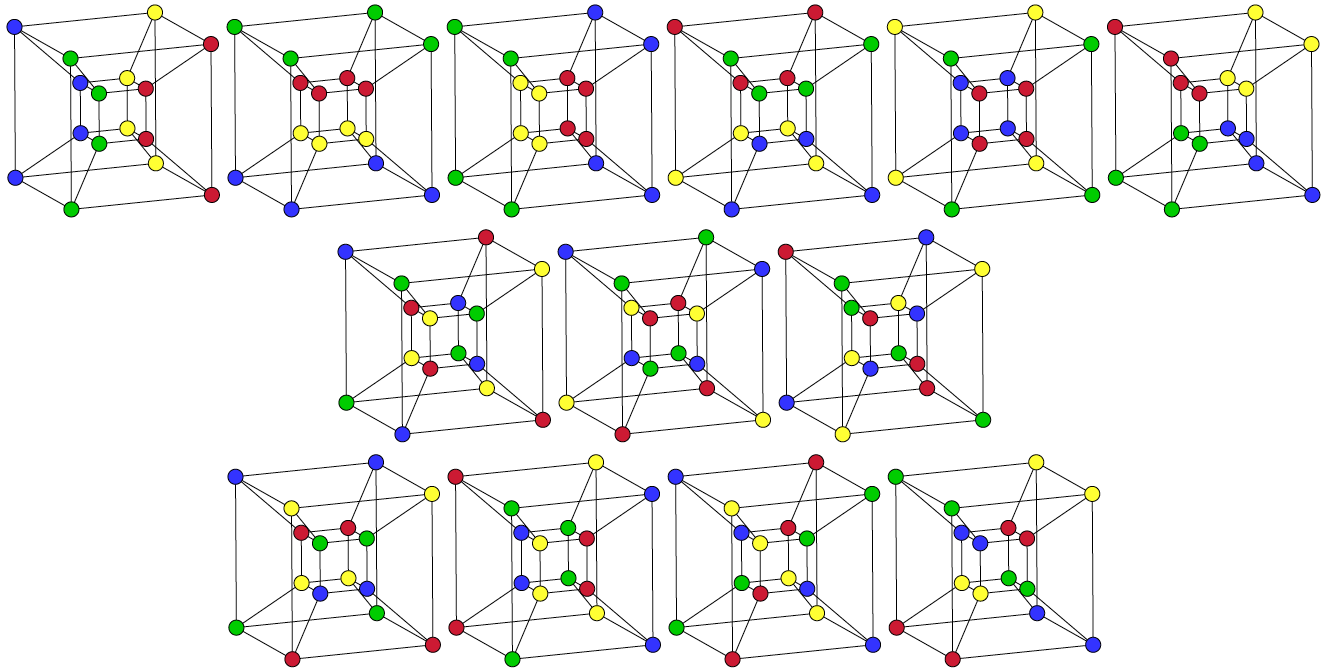}
\end{figure}

\noindent{\bf The symmetries of the hypercube:} It is easy to describe the symmetries of the interval, the square, and the cube. For the hypercube, 
the situation is slightly more complicated due to our loss of geometric intuition (since the geometry we see is not the geometry inherent in dimension 
$4$). However, we can reason by analogy. To create a symmetry, we must proceed as follows.

\begin{itemize}
\item Once a vertex has been chosen, we must decide which vertex it is sent to. Obviously, we have $16$ options for this.

\item Next, we must decide what to do with the four neighbors of the initial vertex, which must be sent to the 
neighbors of the image vertex. Alternatively, we must send the $4$ edges starting from the initial vertex to the $4$ edges leaving the image vertex. 
Obviously, there are $4 \times 3 \times 2 \times 1 = 24$ ways to make this distribution.
\end{itemize}

It is not difficult to convince oneself that once these choices 
are made, the symmetry is completely determined (that is to say, there is only one way to ``extend" it). Consequently, the hypercube admits exactly 
$16 \times 24 = 384$ symmetries.

As a matter of example, below we represent an involution (that is, a symmetry of order 2) of the hypercube. In fact, under an appropriate 
identification, this corresponds to the involution of the $4 \times 4$ square previously introduced.

\begin{figure}[htbp]
    \centering
    \includegraphics[width=0.6\textwidth]{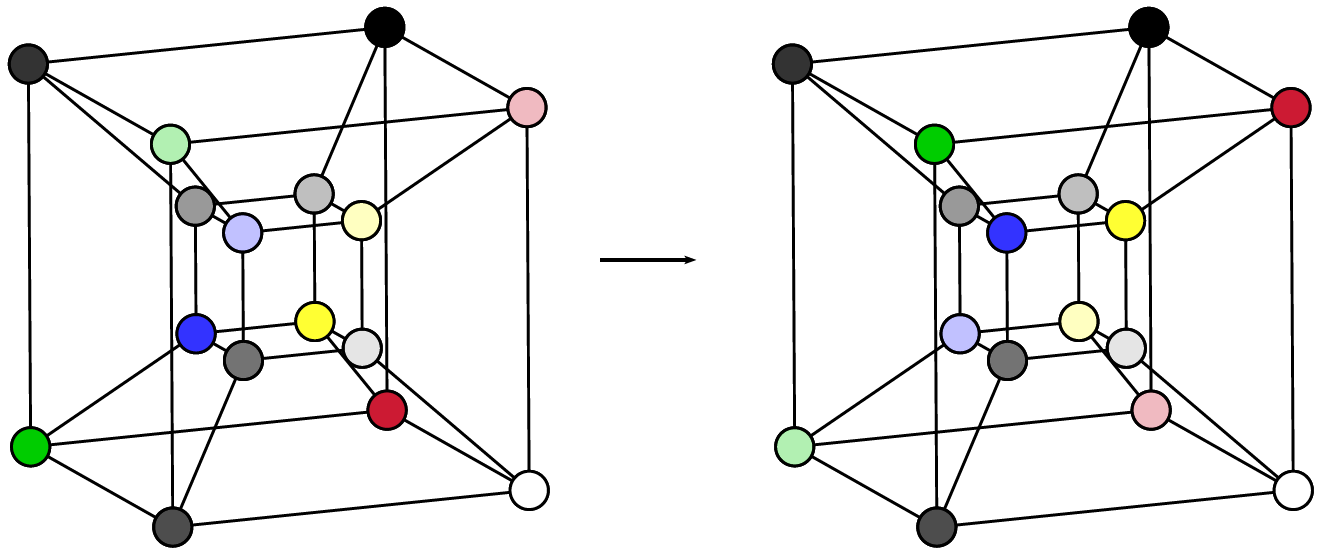}
\end{figure}

\vspace{0.2cm}

\noindent{\bf SUMMARY:}

\begin{itemize}

\item Every symmetry of the hypercube is panmagic when viewed as a permutation 
of the cells of the $4 \times 4$ square (this corresponds to {{the Claim}} above).

\item The hypercube admits exactly 384 symmetries.

\item Applied to the Chautisa Yantra, these symmetries generate $384$ 
distinct panmagic squares with (all different) entries $1, 2, \ldots, 16$.

\item Since there are exactly $384$ panmagic squares with entries $1, 2, \ldots, 16$ 
($\mathrm{N\overline{a}r\overline{a}ya\underaccent{.}na}$'s theorem), there are no more panmagic permutations.
\end{itemize}

We therefore conclude:

\vspace{0.2cm}

\begin{figure}[htbp]
    \centering
    \includegraphics[width=0.7\textwidth]{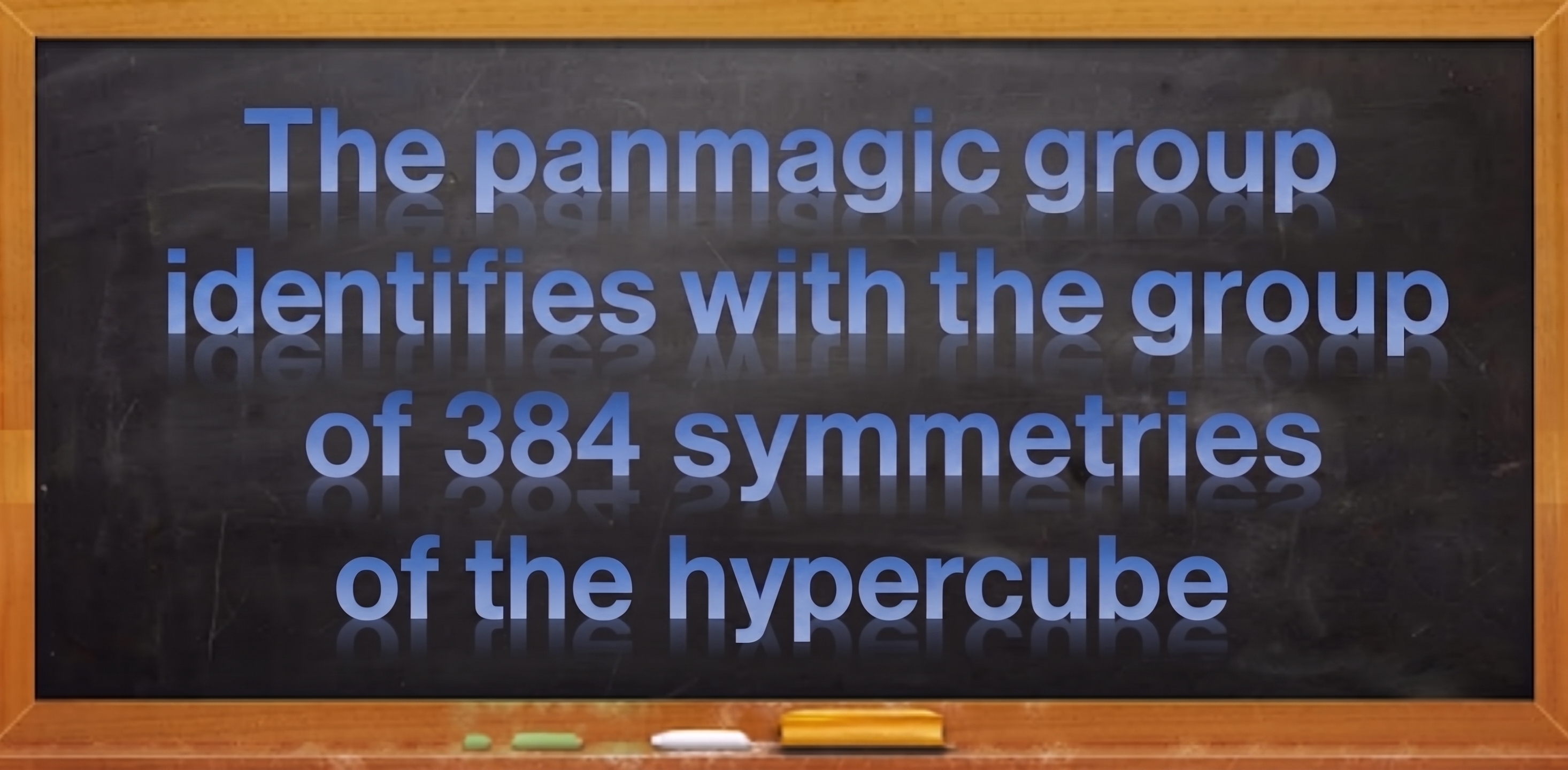}
\end{figure}

\newpage

\noindent{\bf A word on algebraic structure: solvability.} Saying that a group possesses exactly $384$ elements does not bring much, 
because there are $20169$ such non-isomorphic groups.  
However, according to an old 
theorem by W. Burnside, an important property that one draws just from the cardinality is the solvability of the group. 
Without wanting to go into the detail of this notion, we can point out that it is one of the most fundamental in group theory. 
In fact, it dates back to its origins, as it is closely linked to the problem of the solvability of algebraic equations by radicals, 
according to Galois theory.

Now, in our case, the geometric approach allows us to completely describe the group in algebraic terms (and more easily 
corroborate its solvability): it is a semi-direct product
$$S_4 \ltimes \big( \mathbb{Z} / 2 \mathbb{Z}\big)^4.$$
Nothing mysterious here: in the model of the $4 \times 4$ square equipped with the Distance, each factor $\mathbb {Z} _2$ 
corresponds to the action of changing each cell for another in which a specified index changes (geometrically, this corresponds 
to permuting two pairs of adjacent rows or two pairs of adjacent columns). Similarly, the factor $S_4$ is associated with the 
process of ``choosing where the four neighbors of the starting vertex go" described above.

One can spend a pleasant evening unraveling the algebra of this group a bit more by associating its elements with all the 
variations of the Khajuraho square. One should end up with a list of $384$ magic hypercubes~! Here is one to start with:

\begin{figure}[htbp]
    \centering
    \includegraphics[width=0.585\textwidth]{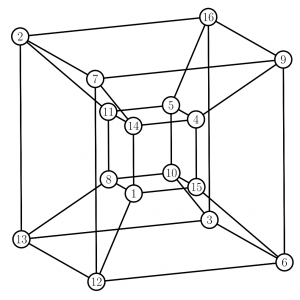}
\end{figure}

\vspace{0.1cm}

\noindent {\bf PAST AND FUTURE OF THE SUBJECT} 

\vspace{0.2cm}

The theorem presented here could be attributed to B. Rosser and R. J. Walker, although they formulate it very differently in their 1938 paper 
\cite{RW38}. In fact, in the comments of their paper, they mention that H. Coxeter had pointed out to them that, very likely, the group they 
were considering is the same group as that of the symmetries of the hypercube (there is no explicit mention of the hypercube in the core 
of the paper, but just a description of the group in question). The proof presented here (as well as that of 
$\mathrm{N\overline{a}r\overline{a}ya\underaccent{.}na}$'s 
theorem) is new, although another argument can be found in the work of W. Müller [Mu97].

\newpage

Magic configurations exist for $n \times n$ squares. For instance, there 
is a famous $4 \times 4$ square with the numbers $1, 2, \ldots, 16$ which is magic but not panmagic, namely, 
the one illustrated by Albrecht Dürer in his engraving {Melencolia I}. This magic square has been extensively discussed in the literature; in 
particular, a short discussion of its symmetries appears at the end of \cite{Hu11}. The conclusion is clear: since this square possesses fewer 
magic properties than that of Khajuraho (exercise: find the combinations that fail...), we end up with fewer symmetries. In our language, 
this translates to the fact that the $4 \times 4$ magic group is ``small": it is reduced to the dihedral group $D_4$ already found for the 
Lo Shu square.\footnote{Why didn't Dürer construct a panmagic square? The reason is simple: he wanted to place the numbers $14$ and 
$15$ next to each other so as to make the number $1514$, which corresponds to the year of the creation of the engraving. However, in 
the proof of $\mathrm{N\overline{a}r\overline{a}ya\underaccent{.}na}$'s theorem, we saw that, in a panmagic square with entries $1, 2, \ldots, 16$, 
the numbers $14$ and $15$ must be neighbors of $1$, and therefore cannot be neighbors of each other.}

\begin{figure}[htbp]
    \centering
    \includegraphics[width=0.475\textwidth]{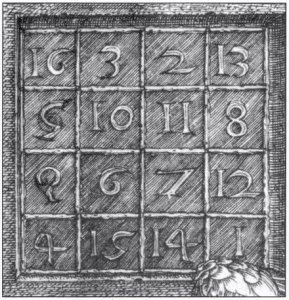}
\end{figure}

Panmagic configurations also exist for other $n \times n$ tables. Of course, magic and panmagic groups 
(and other related groups) can be defined analogously\footnote{In fact, analogous groups can be defined for magic configurations that are not 
necessarily square —such as {\em magic stars}— of numbers.}. Unfortunately, our poor knowledge of higher-dimensional magic squares makes it 
impossible to determine their symmetry groups using the strategy above. For example, the following question seems to me largely open 
(and very interesting):

\vspace{0.35cm}

\noindent{\bf Question.} {\em Are high-dimensional panmagic groups (non-)solvable?}

\vspace{0.35cm}

A more difficult problem is that of the solvability of the image of this group acting on the space of $n \times n$ panmagic 
squares filled with all the numbers $1, 2, \ldots, n^2$ (in case they exist). Starting from $n = 8$, the answer to both 
problems becomes unclear.

\vspace{0.15cm}

To conclude, here is a theorem taken from Müller's paper \cite{Mu97} 
for which I do not know if a more geometric version or proof can be given. 
Note that the presence of $S_5$ in the underlying group implies that it is not solvable.

\vspace{0.35cm}

\noindent{\bf Threorem.} {\em The panmagic group of $5 \times 5$ squares is isomorphic to the semi-direct product $\mathbb{Z}_2 \ltimes (S_5 \times S_5)$}

\vspace{0.25cm}

\noindent{\bf Proof.} Once again, we import from \cite{Na19} the general formula for $5 \times 5$ panmagic squares, whose simplest writing is that of a 
Greco-Latin square:

\begin{figure}[htbp]
    \centering
    \includegraphics[width=0.99\textwidth]{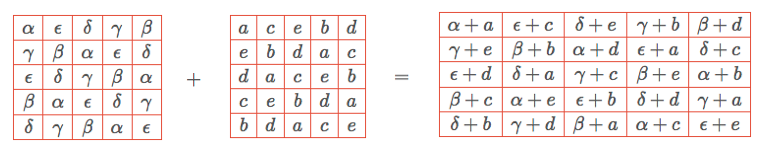}
\end{figure}

\newpage 

The first exercise consists of listing all squares of this type that use the numbers $1, 2, \ldots, 25$. 
For this, we use the following claim, the proof of which is a nice exercise:

\vspace{0.25cm}

\noindent{\bf Claim/exercise.} {\em There are exactly $28800 = 2 \times (5!)^2$ panmagic 
$5 \times 5$ squares that use each of the numbers $1, 2, \ldots, 25$ exactly once. These are obtained 
from the general formula above by setting $\{ \alpha, \beta, \gamma, \delta, \epsilon \} = \{ 0, 5, 10, 15, 20 \}$ 
and $\{ a, b, c, d, e \} = \{1, 2, 3, 4, 5 \}$, or vice-versa.}

\vspace{0.25cm}

Observe now that the group $\mathbb {Z} _2 \ltimes (S_5 \times S_5)$ acts naturally on the set of $5 \times 5$ panmagic 
squares: one factor $S_5$ acts by permuting the entries $\alpha, \beta, \gamma, \delta, \epsilon$, the other factor permutes 
the entries $a, b, c, d, e$, and the factor $\mathbb {Z} _2$ interchanges the Greek and Latin parameters. 
The fact that the panmagic properties are preserved by these permutations is corroborated by a simple 
inspection with the general formula above.

The $5 \times 5$ panmagic group also acts on the set of $5 \times 5$ 
panmagic squares with entries $1, 2, \ldots, 25$. According to the claim/exercise, there are $28800$ such squares. 
Since the number of elements of the group $\mathbb {Z} _2 \ltimes (S_5 \times S_5)$ is also equal to $2 \times (5!)^2 = 28800$, 
we conclude that there can be no other panmagic permutation.

\vspace{0.2cm}

\noindent {\bf To play:} As with the Chautisa Yantra, you can take 
your favorite $5 \times 5$ panmagic square and have fun applying symmetries to it. My favorite square is this 
one:\footnote{This panmagic square has the additional property that entries in positions symmetric with respect 
to the center always sum to $26$, which is twice the central entry (equal to $13$). This square comes from the 
Islamic world, where its central symmetry represents ``the circulation of everything around Allah''.}

\begin{figure}[htbp]
    \centering
    \includegraphics[width=0.29\textwidth]{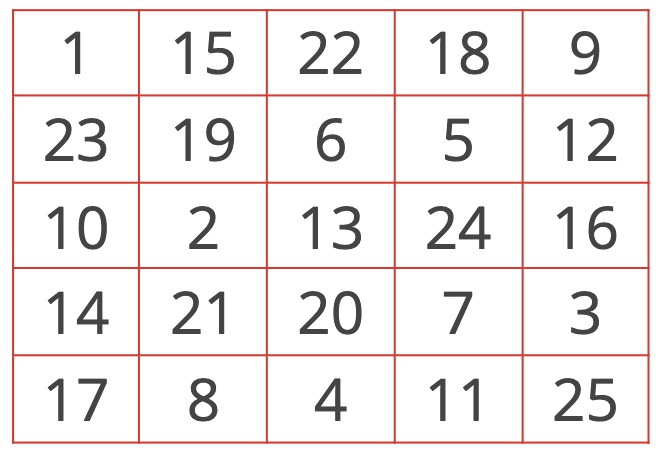}
\end{figure}

\noindent A symmetry that I like very much is this one:

\begin{figure}[htbp]
    \centering
    \includegraphics[width=0.65\textwidth]{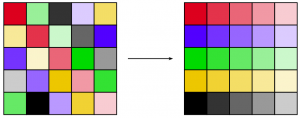}
\end{figure}

\noindent Have fun! There are 28800 panmagic squares waiting for you~!


\begin{footnotesize}

\vspace{0.5cm}

\noindent {\large Andr\'es Navas \\}
\noindent email: andres.navas@usach.cl\\

\noindent Departamento de Matem\'atica y Ciencia de la Computaci\'on\\
\noindent Universidad de Santiago de Chile (USACH)\\ 
\noindent Alameda 3363, Santiago, Chile\\ 

\noindent Institut des Mathématiques\\
Université de Rennes\\
Bat. 22-23 Campus Beaulieu, 35042, France\\

\end{footnotesize}

\end{document}